\documentclass{amsart}
\usepackage{amssymb}
\usepackage{amsmath}
\usepackage{amscd}
\usepackage[left=3cm, right=2.5cm, top=2cm, bottom=2cm]{geometry}
\usepackage{hyperref}

\newtheorem{prop}{Proposition}[section]

\newtheorem{lemma}[prop]{Lemma}
\newtheorem{theorem}[prop]{Theorem}

\newtheorem{notation}[prop]{Notation}
\theoremstyle{definition}

\DeclareFontEncoding{OT2}{}{} 

\begin{document}

\title{On a conjecture of Sun}
\author{Srilakshmi Krishnamoorthy}
\author{Abinash Sarma}

\address{Indian Institute of Science Education and Research Thiruvananthapuram, Maruthamala P.O., Vithura, Thiruvananthapuram-695551, Kerala, India.}
\email{srilakshmi@iisertvm.ac.in}
\address{Indian Institute of Science Education and Research Thiruvananthapuram, Maruthamala P.O., Vithura, Thiruvananthapuram-695551, Kerala, India.}
\email{sarmaabinash15@iisertvm.ac.in}

\begin{abstract}
A number of the form $x(x+1)/2$ where $x$ is an integer is called a triangular number. Suppose, $N(a_1,\cdots,a_k;n)$ and $T(a_1,\cdots,a_k;n)$ denote the number of ways $n$ can be expressed as $\sum_{i=1}^k a_ix_i^2$ and $\sum_{i=1}^k a_i\frac{x_i(x_i+1)}{2}$, respectively. Z.-H. Sun, in \cite{4}, conjectured some relations between $T(a,b,c;n)$ and $N(a,b,c;8n+a+b+c)$. In this paper, we prove these conjectures using theta function identities. Moreover, we add some new triplets $(a,b,c)$ satisfying these conjectures.
\end{abstract}

\maketitle

\section{Introduction}
Let $\mathbb{Z}$, $\mathbb{Z}^+$, and $\mathbb{N}$ denote the set of integers, the set of positive integers, and the set of non-negative integers, respectively. The numbers of the form $x(x+1)/2$ where $x\in\mathbb{Z}$ are called triangular numbers. For $n\in\mathbb{N}$ and $a_1,a_2,\cdots,a_k\in\mathbb{Z}^+$ ($k\geq2$), we set the following notations.
\begin{align*}
    N(a_1,a_2,\cdots,a_k;n):=&\#\left\{(x_1,x_2,\cdots,x_k)\in\mathbb{Z}^k:n=\sum_{i=1}^ka_ix_i^2\right\},\\
    T(a_1,a_2,\cdots,a_k;n):=&\#\left\{(x_1,x_2,\cdots,x_k)\in\mathbb{N}^k:n=\sum_{i=1}^ka_i\frac{x_i(x_i+1)}{2}\right\}.
\end{align*}

Also, define
\begin{align*}
    C(a_1,a_2,\cdots,a_k):=\binom{i_1}{4}+\binom{i_1}{2}+i_1i_3,
\end{align*}
where $i_j$ denotes the number of times $j$ appears in $(a_1,a_2,\cdots,a_k)$. Adiga, Cooper, and Han \cite{1} showed that for $n\in\mathbb{N}$ and for $a_1+a_2+\cdots+a_k\leq7$,
\begin{align*}
    2^kT(a_1,a_2,\cdots,a_k;n)=\frac{2}{2+C(a_1,a_2,\cdots,a_k)}N(a_1,a_2,\cdots,a_k;8n+a_1+\cdots+a_k).
\end{align*}
Baruah, Cooper, and Hirchhorn \cite{2} proved that for $a_1+a_2+\cdots+a_k=8$,
\begin{align*}
    2^kT(a_1,a_2,\cdots,a_k;n)=\frac{2}{2+C(a_1,a_2,\cdots,a_k)}\big(&N(a_1,a_2,\cdots,a_k;8n+8)\\
    &-N(a_1,a_2,\cdots,a_k;2n+2)\big).
\end{align*}
For a detailed literature review, please look at \cite{4}.

Sun \cite{4} derived some relations between $T(a,b,c;n)$ and $N(a,b,c;8n+a+b+c)$. In the last section of this paper, he also stated the following conjectures based on computational evidence for $n\leq150$.

\noindent\textbf{Conjecture 6.1.} \textit{Let $n\in\mathbb{Z}^+$. For $(a,b,c)$ $=$ $(1,1,3)$, $(1,1,4)$, $(1,1,6)$, $(1,1,7)$, $(1,1,15)$, $(1,2,2)$, $(1,2,5)$, $(1,3,3)$, $(1,3,9)$, $(1,5,10)$, $(1,6,9)$, $(1,7,7)$, $(1,7,15)$, $(1,9,15)$, $(1,15,15)$, $(1,15,25)$, $(2,3,3)$ we have}
\begin{align*}
    16T(a,b,c;n)=N(a,b,c;4(8n+a+b+c))-N(a,b,c;8n+a+b+c).
\end{align*}
\noindent\textbf{Conjecture 6.2.} \textit{Let $n\in\mathbb{Z}^+$.\\
(i) For even $n$ and $(a,b,c)$ $=$ $(1,2,15)$, $(1,15,18)$, $(1,15,30)$, $(3,10,45)$ we have}
\begin{align*}
    16T(a,b,c;n)=N(a,b,c;4(8n+a+b+c))-N(a,b,c;8n+a+b+c).
\end{align*}
\textit{(ii) For odd $n$ and $(a,b,c)$ $=$ $(1,6,7)$, $(1,7,42)$, $(2,3,21)$, $(2,9,15)$, $(3,5,6)$, $(3,5,10)$, $(5,21,35)$ we have}
\begin{align*}
    16T(a,b,c;n)=N(a,b,c;4(8n+a+b+c))-N(a,b,c;8n+a+b+c).
\end{align*}
\noindent\textbf{Conjecture 6.3.} \textit{Let $n\in\mathbb{Z}^+$. For $(a,b,c)$ $=$ $(1,3,5)$, $(1,3,7)$, $(1,3,15)$, $(1,3,21)$, $(1,5,15)$, $(1,7,21)$, $(3,5,9)$, $(3,5,15)$, $(3,7,21)$ we have}
\begin{align*}
    16T(a,b,c;n)=3N(a,b,c;8n+a+b+c)-N(a,b,c;4(8n+a+b+c)).
\end{align*}
\noindent\textbf{Conjecture 6.4.} \textit{Let $n\in\mathbb{Z}^+$.\\
(i) For even $n$ and $(a,b,c)$ $=$ $(1,6,15)$, $(1,10,15)$ we have}
\begin{align*}
    16T(a,b,c;n)=3N(a,b,c;8n+a+b+c)-N(a,b,c;4(8n+a+b+c)).
\end{align*}
\textit{(ii) For odd $n$ and $(a,b,c)$ $=$ $(1,2,7)$, $(1,7,14)$, $(2,3,5)$, $(3,5,10)$ we have}
\begin{align*}
    16T(a,b,c;n)=3N(a,b,c;8n+a+b+c)-N(a,b,c;4(8n+a+b+c)).
\end{align*}

\noindent Recently, Xia and Yan \cite{5} proved Conjecture 6.1 for $(1,1,7)$, $(1,1,15)$, $(1,7,7)$, $(1,7,15)$, $(1,9,15)$, $(1,15,15)$, $(1,15,25)$; Conjecture 6.2 (i) for $(1,2,15)$, $(1,15,18)$, $(1,15,30)$; and Conjecture 6.3. For more related works, please look at \cite{22,6,7,8,9,10,11,12,13,14,15,16,17,18,19,20,21}.

The aim of this paper is to prove Conjecture 6.1 for $(1,1,3)$, $(1,1,4)$, $(1,1,6)$, $(1,2,2)$, $(1,3,3)$, $(1,3,9)$, $(2,3,3)$; Conjecture 6.2 (i) for $(3,6,45)$, $(3,10,45)$; Conjecture 6.2 (ii); Conjecture 6.3 for $(3,7,15)$; and Conjecture 6.4. Furthermore, we will prove that $(5,21,35)$, in fact, satisfies Conjecture 6.1. To be precise, the following are the main results of this paper.
\begin{theorem}\label{thm1}
Let $n\in\mathbb{Z}^+$. For $(a,b,c)$ $=$ $(1,1,3)$, $(1,1,4)$, $(1,1,6)$, $(1,2,2)$, $(1,3,3)$, $(1,3,9)$, $(2,3,3)$, $(5,21,35)$ we have
\begin{align*}
    16T(a,b,c;n)=N(a,b,c;4(8n+a+b+c))-N(a,b,c;8n+a+b+c).
\end{align*}
\end{theorem}
\begin{theorem}\label{thm2}
Let $n\in\mathbb{Z}^+$.\\
(i) For even $n$ and $(a,b,c)$ $=$ $(3,6,45)$, $(3,10,45)$ we have
\begin{align*}
    16T(a,b,c;n)=N(a,b,c;4(8n+a+b+c))-N(a,b,c;8n+a+b+c).
\end{align*}
(ii) For odd $n$ and $(a,b,c)$ $=$ $(1,6,7)$, $(1,7,42)$, $(2,3,21)$, $(2,9,15)$, $(3,5,6)$, $(3,5,10)$ we have
\begin{align*}
    16T(a,b,c;n)=N(a,b,c;4(8n+a+b+c))-N(a,b,c;8n+a+b+c).
\end{align*}
\end{theorem}
\begin{theorem}\label{thm3}
Let $n\in\mathbb{Z}^+$. For $(a,b,c)$ $=$ $(3,7,15)$ we have
\begin{align*}
    16T(a,b,c;n)=3N(a,b,c;8n+a+b+c)-N(a,b,c;4(8n+a+b+c)).
\end{align*}
\end{theorem}
\begin{theorem}\label{thm4}
Let $n\in\mathbb{Z}^+$.\\
(i) For even $n$ and $(a,b,c)$ $=$ $(1,6,15)$, $(1,10,15)$ we have
\begin{align*}
    16T(a,b,c;n)=3N(a,b,c;8n+a+b+c)-N(a,b,c;4(8n+a+b+c)).
\end{align*}
(ii) For odd $n$ and $(a,b,c)$ $=$ $(1,2,7)$, $(1,7,14)$, $(2,3,5)$, $(3,5,30)$ we have
\begin{align*}
    16T(a,b,c;n)=3N(a,b,c;8n+a+b+c)-N(a,b,c;4(8n+a+b+c)).
\end{align*}
\end{theorem}

\section{Preliminaries}
Ramanujan's theta functions $\varphi(q)$ and $\psi(q)$ are defined by
\begin{align*}
    \varphi(q)=\sum_{n=-\infty}^\infty q^{n^2},
\end{align*}
and
\begin{align*}
    \psi(q)=\sum_{n=0}^\infty q^{n(n+1)/2}.
\end{align*}
It is easy to see that the generating functions for $N(a_1,a_2,\cdots,a_k;n)$ and $T(a_1,a_2,\cdots,a_k;n)$ are given by
\begin{align}
    \sum_{n=0}^\infty N(a_1,a_2,\cdots,a_k;n)q^n=\varphi(q^{a_1})\varphi(q^{a_2})\cdots\varphi(q^{a_k}),\label{eq:1}\\
    \sum_{n=0}^\infty T(a_1,a_2,\cdots,a_k;n)q^n=\psi(q^{a_1})\psi(q^{a_2})\cdots\psi(q^{a_k}).\label{eq:2}
\end{align}

Let's list the following identities involving $\varphi(q)$ and $\psi(q)$, which will be used in the proofs.

\begin{lemma}\label{lemma2}
The following identities hold:
\begin{align}
    \varphi(q)\varphi(q^{15})&=\varphi(-q^6)\varphi(-q^{10})+2q\psi(q^3)\psi(q^5),\label{eq:11}\\
    \varphi(q^3)\varphi(q^5)&=\varphi(-q^2)\varphi(-q^{30})+2q^2\psi(q)\psi(q^{15}),\label{eq:12}\\
    2\psi(q^6)\psi(q^{10})&=\psi(q)\psi(q^{15})+\psi(-q)\psi(-q^{15}),\label{eq:13}\\
    2q^3\psi(q^2)\psi(q^{30})&=\psi(q^3)\psi(q^5)-\psi(-q^3)\psi(-q^5).\label{eq:10}
\end{align}
\end{lemma}
\noindent These identities from Ramanujan's notebooks appear on page no. 377 in \cite{3} with proof.

\begin{lemma}\label{lemma1}
The following identities hold:
\begin{align}
    \varphi(q)&=\varphi(q^4)+2q\psi(q^8),\label{eq:5}\\
    \varphi(q)^2&=\varphi(q^2)^2+4q\psi(q^4)^2,\label{eq:37}\\
    \psi(q)^2&=\psi(q^2)\varphi(q^4)+2q\psi(q^2)\psi(q^8),\label{eq:38}\\
    \psi(q)\psi(q^3)&=\psi(q^4)\varphi(q^6)+q\varphi(q^2)\psi(q^{12}),\label{eq:19}\\
    \varphi(q)\varphi(q^3)&=\varphi(q^4)\varphi(q^{12})+2q\psi(q^2)\psi(q^6)+4q^4\psi(q^8)\psi(q^{24}),\label{eq:20}\\
    \psi(q)\psi(q^7)&=\psi(q^8)\varphi(q^{28})+q\psi(q^2)\psi(q^{14})+q^6\varphi(q^4)\psi(q^{56}),\label{eq:35}\\
    \psi(q)\psi(q^{15})&=\psi(q^6)\psi(q^{10})+q\varphi(q^{20})\psi(q^{24})+q^3\varphi(q^{12})\psi(q^{40}),\label{eq:3}\\
    \psi(q^3)\psi(q^5)&=\psi(q^8)\varphi(q^{60})+q^3\psi(q^2)\psi(q^{30})+q^{14}\varphi(q^4)\psi(q^{120}).\label{eq:8}
\end{align}
\end{lemma}
\noindent The proof of this lemma is omitted; please look at \cite{5} for the same.

\begin{notation}
Let $h(q)=\sum_{n=0}^\infty c(n)q^n$ be a power series. Denote
\begin{align*}
    \mathcal{E}(h(q)):=&\sum_{n=0}^\infty c(2n)q^{2n},\\
    \mathcal{O}(h(q)):=&\sum_{n=0}^\infty c(2n+1)q^{2n+1}.
\end{align*}
\end{notation}

\section{Proof of Theorem \ref{thm1}}
We will only prove Theorem \ref{thm1} for $(a,b,c)=(5,21,35)$. The other cases can be proved similarly using Lemma \ref{lemma1}.

By (\ref{eq:2}), we have
\begin{align}
    \sum_{n=0}^\infty T(5,21,35;n)q^n=\psi(q^5)\psi(q^{21})\psi(q^{35}).\label{eq:32}
\end{align}
On the other hand, by (\ref{eq:1}), we have
\begin{align}
    \sum_{n=0}^\infty N(5,21,35;n)q^n=\varphi(q^5)\varphi(q^{21})\varphi(q^{35}).\label{eq:33}
\end{align}
In view of (\ref{eq:5}), we get
\begin{align*}
    \sum_{n=0}^\infty &N(5,21,35;n)q^n\\
    &=(\varphi(q^{20})+2q^5\psi(q^{40}))(\varphi(q^{84})+2q^{21}\psi(q^{168}))(\varphi(q^{140})+2q^{35}\psi(q^{280})).
\end{align*}
Extracting the terms involving $q^{4n}$ and substituting $q$ by $q^{1/4}$ yields after rearrangement
\begin{align}
    \sum_{n=0}^\infty N(5,21,35;4n)q^n&=\varphi(q^5)\varphi(q^{21})\varphi(q^{35})+4q^{10}\psi(q^{10})\varphi(q^{21})\psi(q^{70})\nonumber\\
    &\quad+4q^{14}\varphi(q^5)\psi(q^{42})\psi(q^{70}).\label{eq:34}
\end{align}
We define
\begin{align*}
    r(n):=N(5,21,35;4n)-N(5,21,35;n).
\end{align*}
By (\ref{eq:33}) and (\ref{eq:34}), we have
\begin{align*}
    \sum_{n=0}^\infty r(n)q^n=4q^{10}\psi(q^{10})\varphi(q^{21})\psi(q^{70})+4q^{14}\varphi(q^5)\psi(q^{42})\psi(q^{70}).
\end{align*}
In view of (\ref{eq:5}), we get
\begin{align*}
    \sum_{n=0}^\infty &r(2n+1)q^{2n+1}\\
    &=4q^{10}\psi(q^{10})\psi(q^{70})\mathcal{O}(\varphi(q^{21}))+4q^{14}\psi(q^{42})\psi(q^{70})\mathcal{O}(\varphi(q^5))\\
    &=4q^{10}\psi(q^{10})\psi(q^{70})(2q^{21}\psi(q^{168}))+4q^{14}\psi(q^{42})\psi(q^{70})(2q^5\psi(q^{40})).
\end{align*}
Dividing by $q$ and substituting $q$ by $q^{1/2}$ yields after rearrangement
\begin{align*}
    \sum_{n=0}^\infty r(2n+1)q^n=8q^9\psi(q^{20})\psi(q^{21})\psi(q^{35})+8q^{15}\psi(q^5)\psi(q^{35})\psi(q^{84}).
\end{align*}
In view of (\ref{eq:35}) and (\ref{eq:8}), we get
\begin{align*}
    \sum_{n=0}^\infty r(4n+1)q^{2n}&=8q^9\psi(q^{20})\mathcal{O}(\psi(q^{21})\psi(q^{35}))+8q^{15}\psi(q^{84})\mathcal{O}(\psi(q^5)\psi(q^{35}))\\
    &=8q^9\psi(q^{20})(q^{21}\psi(q^{14})\psi(q^{210}))+8q^{15}\psi(q^{84})(q^5\psi(q^{10})\psi(q^{70})).
\end{align*}
Substituting $q$ by $q^{1/2}$ yields after rearrangement
\begin{align*}
    \sum_{n=0}^\infty r(4n+1)q^n=8q^{10}\psi(q^5)\psi(q^{35})\psi(q^{42})+8q^{15}\psi(q^7)\psi(q^{10})\psi(q^{105}).
\end{align*}
In view of (\ref{eq:35}) and (\ref{eq:3}), we get
\begin{align*}
    \sum_{n=0}^\infty r(8n+5)q^{2n+1}&=8q^{10}\psi(q^{42})\mathcal{O}(\psi(q^5)\psi(q^{35}))+8q^{15}\psi(q^{10})\mathcal{E}(\psi(q^7)\psi(q^{105}))\\
    &=8q^{10}\psi(q^{42})(q^5\psi(q^{10})\psi(q^{70}))+8q^{15}\psi(q^{10})(\psi(q^{42})\psi(q^{70})).
\end{align*}
Dividing by $q$ and substituting $q$ by $q^{1/2}$ yields after rearrangement
\begin{align}
    \sum_{n=0}^\infty r(8n+5)q^n=16q^7\psi(q^5)\psi(q^{21})\psi(q^{35}).\label{eq:36}
\end{align}
By (\ref{eq:32}) and (\ref{eq:36}), we get
\begin{align*}
    16\sum_{n=0}^\infty T(5,21,35;n)q^{n+7}=\sum_{n=0}^\infty r(8n+5)q^n.
\end{align*}
Comparing the coefficients of $q^{n+7}$, we can conclude that
\begin{align*}
    16T(5,21,35;n)=N(5,21,35;4(8n+61))-N(5,21,35;8n+61).
\end{align*}

\section{Proof of Theorem \ref{thm2}}
\textbf{Part (i)}: We will only prove Theorem \ref{thm2} (i) for $(a,b,c)=(3,10,45)$. The other cases can be proved similarly using Lemma \ref{lemma1}.

By (\ref{eq:2}), we have
\begin{align*}
    \sum_{n=0}^\infty T(3,10,45;n)q^n=\psi(q^3)\psi(q^{10})\psi(q^{45}).
\end{align*}
In view of (\ref{eq:3}), we get
\begin{align*}
    \sum_{n=0}^\infty T(3,10,45;2n)q^{2n}&=\psi(q^{10})\mathcal{E}(\psi(q^3)\psi(q^{45}))\\
    &=\psi(q^{10})(\psi(q^{18})\psi(q^{30})).
\end{align*}
Substituting $q$ by $q^{1/2}$ yields after rearrangement
\begin{align}
    \sum_{n=0}^\infty T(3,10,45;2n)q^n=\psi(q^5)\psi(q^9)\psi(q^{15}).\label{eq:28}
\end{align}
On the other hand, by (\ref{eq:1}), we have
\begin{align}
    \sum_{n=0}^\infty N(3,10,45;n)q^n=\varphi(q^3)\varphi(q^{10})\varphi(q^{45}).\label{eq:29}
\end{align}
In view of (\ref{eq:5}), we get
\begin{align*}
    \sum_{n=0}^\infty &N(3,10,45;n)q^n\\
    &=(\varphi(q^{12})+2q^3\psi(q^{24}))(\varphi(q^{40})+2q^{10}\psi(q^{80}))(\varphi(q^{180})+2q^{45}\psi(q^{360})).
\end{align*}
Extracting the terms involving $q^{4n}$ and substituting $q$ by $q^{1/4}$ yields after rearrangement
\begin{align}
    \sum_{n=0}^\infty N(3,10,45;4n)q^n=\varphi(q^3)\varphi(q^{10})\varphi(q^{45})+4q^{12}\psi(q^6)\varphi(q^{10})\psi(q^{90}).\label{eq:30}
\end{align}
We define
\begin{align*}
    s_1(n)=N(3,10,45;4n)-N(3,10,45;n).
\end{align*}
By (\ref{eq:29}) and (\ref{eq:30}), we have
\begin{align*}
    \sum_{n=0}^\infty s_1(n)q^n=4q^{12}\psi(q^6)\varphi(q^{10})\psi(q^{90}).
\end{align*}
Clearly,
\begin{align*}
    \sum_{n=0}^\infty s_1(2n)q^{2n}=4q^{12}\psi(q^6)\varphi(q^{10})\psi(q^{90}).
\end{align*}
Substituting $q$ by $q^{1/2}$ yields after rearrangement
\begin{align*}
    \sum_{n=0}^\infty s_1(2n)q^n=4q^6\psi(q^3)\varphi(q^5)\psi(q^{45}).
\end{align*}
In view of (\ref{eq:5}) and (\ref{eq:3}), we get
\begin{align*}
    \sum_{n=0}^\infty &s_1(4n+2)q^{2n+1}\\
    &=4q^6\mathcal{O}(\varphi(q^5))\mathcal{E}(\psi(q^3)\psi(q^{45}))+4q^6\mathcal{E}(\varphi(q^5))\mathcal{O}(\psi(q^3)\psi(q^{45}))\\
    &=4q^6(2q^5\psi(q^{40}))(\psi(q^{18})\psi(q^{30}))+4q^6(\varphi(q^{20}))(q^3\varphi(q^{60})\psi(q^{72})\\
    &\quad+q^9\varphi(q^{36})\psi(q^{120})).
\end{align*}
Dividing by $q$ and substituting $q$ by $q^{1/2}$ yields after rearrangement
\begin{align*}
    \sum_{n=0}^\infty s_1(4n+2)q^n&=4q^4\varphi(q^{10})\varphi(q^{30})\psi(q^{36})+8q^5\psi(q^9)\psi(q^{15})\psi(q^{20})\\
    &\quad+4q^7\varphi(q^{10})\varphi(q^{18})\psi(q^{60}).
\end{align*}
In view of (\ref{eq:8}), we get
\begin{align*}
    \sum_{n=0}^\infty s_1(8n+2)q^{2n}&=4q^4\varphi(q^{10})\varphi(q^{30})\psi(q^{36})+8q^5\psi(q^{20})\mathcal{O}(\psi(q^9)\psi(q^{15}))\\
    &=4q^4\varphi(q^{10})\varphi(q^{30})\psi(q^{36})+8q^5\psi(q^{20})(q^9\psi(q^6)\psi(q^{90})).
\end{align*}
Substituting $q$ by $q^{1/2}$ yields after rearrangement
\begin{align*}
    \sum_{n=0}^\infty s_1(8n+2)q^n=4q^2\varphi(q^5)\varphi(q^{15})\psi(q^{18})+8q^7\psi(q^3)\psi(q^{10})\psi(q^{45}).
\end{align*}
In view of (\ref{eq:20}) and (\ref{eq:3}), we get
\begin{align*}
    \sum_{n=0}^\infty s_1(16n+10)q^{2n+1}&=4q^2\psi(q^{18})\mathcal{O}(\varphi(q^5)\varphi(q^{15}))+8q^7\psi(q^{10})\mathcal{E}(\psi(q^3)\psi(q^{45}))\\
    &=4q^2\psi(q^{18})(2q^5\psi(q^{10})\psi(q^{30}))+8q^7\psi(q^{10})(\psi(q^{18})\psi(q^{30})).
\end{align*}
Dividing by $q$ and substituting $q$ by $q^{1/2}$ yields after rearrangement
\begin{align}
    \sum_{n=0}^\infty s_1(16n+10)q^n=16q^3\psi(q^5)\psi(q^9)\psi(q^{15}).\label{eq:31}
\end{align}
By (\ref{eq:28}) and (\ref{eq:31}), we get
\begin{align*}
    16\sum_{n=0}^\infty T(3,10,45;2n)q^{n+3}=\sum_{n=0}^\infty s_1(16n+10)q^n.
\end{align*}
Comparing the coefficients of $q^{n+3}$, we can conclude that
\begin{align*}
    16T(3,10,45;2n)=N(3,10,45;4(16n+58))-N(3,10,45;16n+58).
\end{align*}

\noindent\textbf{Part (ii)}: We will only prove Theorem \ref{thm2} (ii) for $(a,b,c)=(1,6,7)$. The other cases can be proved similarly using Lemma \ref{lemma1}.

By (\ref{eq:2}), we have
\begin{align*}
    \sum_{n=0}^\infty T(1,6,7;n)q^n=\psi(q)\psi(q^6)\psi(q^7).
\end{align*}
In view of (\ref{eq:35}), we get
\begin{align*}
    \sum_{n=0}^\infty T(1,6,7;2n+1)q^{2n+1}&=\psi(q^6)\mathcal{O}(\psi(q)\psi(q^7))\\
    &=\psi(q^6)(q\psi(q^2)\psi(q^{14})).
\end{align*}
Dividing by $q$ and substituting $q$ by $q^{1/2}$ yields after rearrangement
\begin{align}
    \sum_{n=0}^\infty T(1,6,7;2n+1)q^n=\psi(q)\psi(q^3)\psi(q^7).\label{eq:16}
\end{align}
On the other hand, by (\ref{eq:1}), we have
\begin{align}
    \sum_{n=0}^\infty N(1,6,7;n)q^n=\varphi(q)\varphi(q^6)\varphi(q^7).\label{eq:17}
\end{align}
In view of (\ref{eq:5}), we get
\begin{align*}
    \sum_{n=0}^\infty N(1,6,7;n)q^n=(\varphi(q^4)+2q\psi(q^8))(\varphi(q^{24})+2q^6\psi(q^{48}))(\varphi(q^{28})+2q^7\psi(q^{56})).
\end{align*}
Extracting the terms involving $q^{4n}$ and substituting $q$ by $q^{1/4}$ yields after rearrangement
\begin{align}
    \sum_{n=0}^\infty N(1,6,7;4n)q^n=\varphi(q)\varphi(q^6)\varphi(q^7)+4q^2\psi(q^2)\varphi(q^6)\psi(q^{14}).\label{eq:18}
\end{align}
We define
\begin{align*}
    s_2(n):=N(1,6,7;4n)-N(1,6,7;n).
\end{align*}
By (\ref{eq:17}) and (\ref{eq:18}), we have
\begin{align*}
    \sum_{n=0}^\infty s_2(n)q^n=4q^2\psi(q^2)\varphi(q^6)\psi(q^{14}).
\end{align*}
Clearly,
\begin{align*}
    \sum_{n=0}^\infty s_2(2n)q^{2n}=4q^2\psi(q^2)\varphi(q^6)\psi(q^{14}).
\end{align*}
Dividing by $q$ and substituting $q$ by $q^{1/2}$ yields after rearrangement
\begin{align*}
    \sum_{n=0}^\infty s_2(2n)q^n=4q\psi(q)\varphi(q^3)\psi(q^7).
\end{align*}
In view of (\ref{eq:5}) and (\ref{eq:35}), we get
\begin{align*}
    \sum_{n=0}^\infty &s_2(4n+2)q^{2n+1}\\
    &=4q\mathcal{E}(\varphi(q^3))\mathcal{E}(\psi(q)\psi(q^7))+4q\mathcal{O}(\varphi(q^3))\mathcal{O}(\psi(q)\psi(q^7))\\
    &=4q(\varphi(q^{12}))(\psi(q^8)\varphi(q^{28})+q^6\varphi(q^4)\psi(q^{56}))+4q(2q^3\psi(q^{24}))(q\psi(q^2)\psi(q^{14})).
\end{align*}
Dividing by $q$ and substituting $q$ by $q^{1/2}$ yields after rearrangement
\begin{align*}
    \sum_{n=0}^\infty s_2(4n+2)q^n&=4\psi(q^4)\varphi(q^6)\varphi(q^{14})+8q^2\psi(q)\psi(q^7)\psi(q^{12})\\
    &\quad+4q^3\varphi(q^2)\varphi(q^6)\psi(q^{28}).
\end{align*}
In view of (\ref{eq:35}), we get
\begin{align*}
    \sum_{n=0}^\infty s_2(8n+6)q^{2n+1}&=8q^2\psi(q^{12})\mathcal{O}(\psi(q)\psi(q^7))+4q^3\varphi(q^2)\varphi(q^6)\psi(q^{28})\\
    &=8q^2\psi(q^{12})(q\psi(q^2)\psi(q^{14}))+4q^3\varphi(q^2)\varphi(q^6)\psi(q^{28}).
\end{align*}
Dividing by $q$ and substituting $q$ by $q^{1/2}$ yields after rearrangement
\begin{align*}
    \sum_{n=0}^\infty s_2(8n+6)q^n=4q\varphi(q)\varphi(q^3)\psi(q^{14})+8q\psi(q)\psi(q^6)\psi(q^7).
\end{align*}
In view of (\ref{eq:20}) and (\ref{eq:35}), we get
\begin{align*}
    \sum_{n=0}^\infty s_2(16n+6)q^{2n}&=4q\psi(q^{14})\mathcal{O}(\varphi(q)\varphi(q^3))+8q\psi(q^6)\mathcal{O}(\psi(q)\psi(q^7))\\
    &=4q\psi(q^{14})(2q\psi(q^2)\psi(q^6))+8q\psi(q^6)(q\psi(q^2)\psi(q^{14})).
\end{align*}
Substituting $q$ by $q^{1/2}$ yields after rearrangement
\begin{align}
    \sum_{n=0}^\infty s_2(16n+6)q^n=16q\psi(q)\psi(q^3)\psi(q^7).\label{eq:21}
\end{align}
By (\ref{eq:16}) and (\ref{eq:21}), we get
\begin{align*}
    16\sum_{n=0}^\infty T(1,6,7;n)q^{n+1}=\sum_{n=0}^\infty s_2(16n+6)q^n.
\end{align*}
Comparing the coefficients of $q^{n+1}$, we can conclude that
\begin{align*}
    16T(1,6,7;n)=N(1,6,7;4(16n+22))-N(1,6,7;16n+22).
\end{align*}

\section{Proof of Theorem \ref{thm3}}
By (\ref{eq:2}), we have
\begin{align}
    \sum_{n=0}^\infty T(3,7,15;n)q^n=\psi(q^3)\psi(q^7)\psi(q^{15}).\label{eq:39}
\end{align}
On the other hand, by (\ref{eq:1}), we have
\begin{align}
    \sum_{n=0}^\infty N(3,7,15;n)q^n=\varphi(q^3)\varphi(q^7)\varphi(q^{15}).\label{eq:40}
\end{align}
In view of (\ref{eq:5}), we get
\begin{align}
    \sum_{n=0}^\infty &N(3,7,15;n)q^n\nonumber\\
    &=(\varphi(q^{12})+2q^3\psi(q^{24}))(\varphi(q^{28})+2q^7\psi(q^{56}))(\varphi(q^{60})+2q^{15}\psi(q^{120})).\label{eq:42}
\end{align}
Extracting the terms involving $q^{4n}$ and substituting $q$ by $q^{1//4}$ yields after rearrangement
\begin{align}
    \sum_{n=0}^\infty N(3,7,15;4n)q^n=\varphi(q^3)\varphi(q^7)\varphi(q^{15}).\label{eq:41}
\end{align}
We define
\begin{align*}
    t(n):=3N(3,7,15;n)-N(3,7,15;4n).
\end{align*}
By (\ref{eq:40}) and (\ref{eq:41}), we have
\begin{align*}
    \sum_{n=0}^\infty t(n)q^n=2\varphi(q^3)\varphi(q^7)\varphi(q^{15}).
\end{align*}
In view of (\ref{eq:42}), we get
\begin{align*}
    &\sum_{n=0}^\infty t(2n+1)q^{2n+1}\\
    &=2\mathcal{O}(\varphi(q^3)\varphi(q^7)\varphi(q^{15}))\\
    &=2(2q^3\psi(q^{24})\varphi(q^{28})\varphi(q^{60})+2q^7\varphi(q^{12})\psi(q^{56})\varphi(q^{60})+2q^{15}\varphi(q^{12})\varphi(q^{28})\psi(q^{120})\\
    &\quad+8q^{25}\psi(q^{24})\psi(q^{56})\psi(q^{120})).
\end{align*}
Dividing by $q$ and substituting $q$ by $q^{1/2}$ yields after rearrangement
\begin{align*}
    \sum_{n=0}^\infty t(2n+1)q^n&=4q\psi(q^{12})\varphi(q^{14})\varphi(q^{30})+4q^3\varphi(q^6)\psi(q^{28})\varphi(q^{30})\\
    &\quad+4q^7\varphi(q^6)\varphi(q^{14})\psi(q^{60})+16q^{12}\psi(q^{12})\psi(q^{28})\psi(q^{60}).
\end{align*}
Clearly,
\begin{align*}
    \sum_{n=0}^\infty t(4n+1)q^{2n}=16q^{12}\psi(q^{12})\psi(q^{28})\psi(q^{60}).
\end{align*}
Substituting $q$ by $q^{1/2}$ yields after rearrangement
\begin{align*}
    \sum_{n=0}^\infty t(4n+1)q^n=16q^6\psi(q^6)\psi(q^{14})\psi(q^{30}).
\end{align*}
Clearly,
\begin{align*}
    \sum_{n=0}^\infty t(8n+1)q^{2n}=16q^6\psi(q^6)\psi(q^{14})\psi(q^{30}).
\end{align*}
Substituting $q$ by $q^{1/2}$ yields after rearrangement
\begin{align}
    \sum_{n=0}^\infty t(8n+1)q^n=16q^3\psi(q^3)\psi(q^7)\psi(q^{15}).\label{eq:43}
\end{align}
By (\ref{eq:39}) and (\ref{eq:43}), we have
\begin{align*}
    16\sum_{n=0}^\infty T(3,7,15;n)q^{n+3}=\sum_{n=0}^\infty t(8n+1)q^n.
\end{align*}
Comparing the coefficients of $q^{n+3}$, we can conclude that
\begin{align*}
    16T(3,7,15;n)=3N(3,7,15;n)-N(3,7,15;4n).
\end{align*}

\section{Proof of Theorem \ref{thm4}}
\textbf{Part (i)}: We will only prove Theorem \ref{thm4} (i) for $(a,b,c)$ $=$ $(1,6,15)$. The other cases can be proved similarly using Lemma \ref{lemma2} and Lemma \ref{lemma1}.

By (\ref{eq:2}), we have
\begin{align*}
    \sum_{n=0}^\infty T(1,6,15;n)q^n=\psi(q)\psi(q^6)\psi(q^{15}).
\end{align*}
In view of (\ref{eq:3}), we get
\begin{align*}
    \sum_{n=0}^\infty T(1,6,15;2n)q^{2n}&=\psi(q^6)\mathcal{E}(\psi(q)\psi(q^{15}))\\
    &=\psi(q^6)\psi(q^6)\psi(q^{10}).
\end{align*}
Substituting $q$ by $q^{1/2}$ yields after rearrangement
\begin{align}
    \sum_{n=0}^\infty T(1,6,15;2n)q^n=\psi(q^3)\psi(q^3)\psi(q^5).\label{eq:4}
\end{align}
On the other hand, by (\ref{eq:1}), we have
\begin{align}
    \sum_{n=0}^\infty N(1,6,15;n)q^n=\varphi(q)\varphi(q^6)\varphi(q^{15}).\label{eq:6}
\end{align}
In view of (\ref{eq:5}), we get
\begin{align*}
    \sum_{n=0}^\infty N(1,6,15;n)q^n=(\varphi(q^4)+2q\psi(q^8))(\varphi(q^{24})+2q^6\psi(q^{48}))(\varphi(q^{60})+2q^{15}\psi(q^{120})).
\end{align*}
Extracting the terms involving $q^{4n}$ and substituting $q$ by $q^{1/4}$ yields after rearrangement
\begin{align}
    \sum_{n=0}^\infty N(1,6,15;4n)q^n=\varphi(q)\varphi(q^6)\varphi(q^{15})+4q^4\psi(q^2)\varphi(q^6)\psi(q^{30}).\label{eq:7}
\end{align}
We define
\begin{align*}
    u_1(n):=3N(1,6,15;n)-N(1,6,15;4n).
\end{align*}
By (\ref{eq:6}) and (\ref{eq:7}), we have
\begin{align*}
    \sum_{n=0}^\infty u_1(n)q^n=2\varphi(q)\varphi(q^6)\varphi(q^{15})-4q^4\psi(q^2)\varphi(q^6)\psi(q^{30}).
\end{align*}
In view of (\ref{eq:5}), we get
\begin{align*}
    \sum_{n=0}^\infty u_1(2n)q^{2n}&=2\varphi(q^6)\mathcal{E}(\varphi(q)\varphi(q^{15}))-4q^4\psi(q^2)\varphi(q^6)\psi(q^{30})\\
    &=2\varphi(q^6)(\varphi(q^4)\varphi(q^{60})+4q^{16}\psi(q^8)\psi(q^{120}))-4q^4\psi(q^2)\varphi(q^6)\psi(q^{30}).
\end{align*}
Substituting $q$ by $q^{1/2}$ yields after rearrangement
\begin{align*}
    \sum_{n=0}^\infty u_1(2n)q^n=2\varphi(q^2)\varphi(q^3)\varphi(q^{30})-4q^2\psi(q)\varphi(q^3)\psi(q^{15})+8q^8\varphi(q^3)\psi(q^4)\psi(q^{60}).
\end{align*}
In view of (\ref{eq:5}) and (\ref{eq:3}), we get
\begin{align*}
    \sum_{n=0}^\infty &u_1(4n+2)q^{2n+1}\\
    &=2\varphi(q^2)\varphi(q^{30})\mathcal{O}(\varphi(q^3))-4q^2\mathcal{O}(\varphi(q^3))\mathcal{E}(\psi(q)\psi(q^{15}))\\
    &\quad-4q^2\mathcal{E}(\varphi(q^3))\mathcal{O}(\psi(q)\psi(q^{15}))+8q^8\psi(q^4)\psi(q^{60})\mathcal{O}(\varphi(q^3))\\
    &=2\varphi(q^2)\varphi(q^{30})(2q^3\psi(q^{24}))-4q^2(2q^3\psi(q^{24}))(\psi(q^6)\psi(q^{10}))\\
    &\quad-4q^2(\varphi(q^{12}))(q\varphi(q^{20})\psi(q^{24})+q^3\varphi(q^{12})\psi(q^{40}))+8q^8\psi(q^4)\psi(q^{60})(2q^3\psi(q^{24})).
\end{align*}
Dividing by $q$ and substituting $q$ by $q^{1/2}$ yields after rearrangement
\begin{align*}
    \sum_{n=0}^\infty &u_1(4n+2)q^n\\
    &=4q\varphi(q)\psi(q^2)\varphi(q^{15})-4q\varphi(q^6)\varphi(q^{10})\psi(q^{12})-8q^2\psi(q^3)\psi(q^5)\psi(q^{12})\\
    &\quad-4q^2\varphi(q^6)\varphi(q^6)\psi(q^{20})+16q^5\psi(q^2)\psi(q^{12})\psi(q^{30}).
\end{align*}
In view of (\ref{eq:5}) and (\ref{eq:8}), we get
\begin{align*}
    \sum_{n=0}^\infty &u_1(8n+6)q^{2n+1}\\
    &=4q\psi(q^2)\mathcal{E}(\varphi(q)\varphi(q^{15}))-4q\varphi(q^6)\varphi(q^{10})\psi(q^{12})-8q^2\psi(q^{12})\mathcal{O}(\psi(q^3)\psi(q^5))\\
    &\quad+16q^5\psi(q^2)\psi(q^{12})\psi(q^{30})\\
    &=4q\psi(q^2)(\varphi(q^4)\varphi(q^{60})+4q^{16}\psi(q^8)\psi(q^{120}))-4q\varphi(q^6)\varphi(q^{10})\psi(q^{12})\\
    &\quad-8q^2\psi(q^{12})(q^3\psi(q^2)\psi(q^{30}))+16q^5\psi(q^2)\psi(q^{12})\psi(q^{30}).
\end{align*}
Substituting $q$ by $q^{1/2}$ yields after rearrangement
\begin{align*}
    \sum_{n=0}^\infty u_1(8n+6)q^n&=4\varphi(q^2)\psi(q^6)\varphi(q^{30})-4\varphi(q^3)\varphi(q^5)\psi(q^6)+8q^2\psi(q)\psi(q^6)\psi(q^{15})\\
    &\quad+16q^8\psi(q^4)\psi(q^6)\psi(q^{60}).
\end{align*}
In view of (\ref{eq:5}) and (\ref{eq:3}), we get
\begin{align*}
    \sum_{n=0}^\infty &u_1(16n+6)q^{2n}\\
    &=4\varphi(q^2)\psi(q^6)\varphi(q^{30})-4\psi(q^6)\mathcal{E}(\varphi(q^3)\varphi(q^5))+8q^2\psi(q^6)\mathcal{E}(\psi(q)\psi(q^{15}))\\
    &\quad+16q^8\psi(q^4)\psi(q^6)\psi(q^{60})\\
    &=4\varphi(q^2)\psi(q^6)\varphi(q^{30})-4\psi(q^6)(\varphi(q^{12})\varphi(q^{20})+4q^8\psi(q^{24})\psi(q^{40}))\\
    &\quad+8q^2\psi(q^6)(\psi(q^6)\psi(q^{10}))+16q^8\psi(q^4)\psi(q^6)\psi(q^{60}).
\end{align*}
Substituting $q$ by $q^{1/2}$ yields after rearrangement
\begin{align}
\begin{aligned}
    \sum_{n=0}^\infty &u_1(16n+6)q^n\\
    &=8q\psi(q^3)\psi(q^3)\psi(q^5)\\
    &\quad+4\psi(q^3)\Big(\varphi(q)\varphi(q^{15})-\varphi(q^6)\varphi(q^{10})+4q^4\psi(q^2)\psi(q^{30})-4q^4\psi(q^{12})\psi(q^{20})\Big).\label{eq:9}
\end{aligned}
\end{align}

In view of (\ref{eq:13}) and (\ref{eq:5}), we get
\begin{align*}
    &\varphi(q)\varphi(q^{15})-\varphi(q^6)\varphi(q^{10})+4q^4\psi(q^2)\psi(q^{30})-4q^4\psi(q^{12})\psi(q^{20})\\
    &=\Big(\varphi(q^4)\varphi(q^{60})+2q\psi(q^8)\varphi(q^{60})+2q^{15}\varphi(q^4)\psi(q^{120})+4q^{16}\psi(q^8)\psi(q^{120})\Big)\\
    &\quad-\varphi(q^6)\varphi(q^{10})+4q^4\psi(q^2)\psi(q^{30})-\Big(2q^4\psi(q^2)\psi(q^{30})+2q^4\psi(-q^2)\psi(-q^{30})\Big)\\
    &=2q\Big(\psi(q^8)\varphi(q^{60})+q^3\psi(q^2)\psi(q^{30})+q^{14}\varphi(q^4)\psi(q^{120})\Big)+\varphi(q^4)\varphi(q^{60})\\
    &\quad-\varphi(q^6)\varphi(q^{10})-2q^4\psi(-q^2)\psi(-q^{30})+4q^{16}\psi(q^8)\psi(q^{120}).
\end{align*}
Now by (\ref{eq:8}), we have
\begin{align}
\begin{aligned}[b]
    &\varphi(q)\varphi(q^{15})-\varphi(q^6)\varphi(q^{10})+4q^4\psi(q^2)\psi(q^{30})-4q^4\psi(q^{12})\psi(q^{20})\\
    &=2q\psi(q^3)\psi(q^5)+\varphi(q^4)\varphi(q^{60})-\varphi(q^6)\varphi(q^{10})-2q^4\psi(-q^2)\psi(-q^{30})\\
    &\quad+4q^{16}\psi(q^8)\psi(q^{120}).\label{eq:15}
\end{aligned}
\end{align}
Implementing (\ref{eq:11}), (\ref{eq:12}), (\ref{eq:13}), and (\ref{eq:10}) appropriately, we get
\begin{align*}
    \varphi(q^4)\varphi(q^{60})&=\varphi(-q^{24})\varphi(-q^{40})+2q^4\psi(q^{12})\psi(q^{20}),\\
    \varphi(q^6)\varphi(q^{10})&=\varphi(-q^4)\varphi(-q^{60})+2q^4\psi(q^2)\psi(q^{30})\\
    &=\varphi(-q^{24})\varphi(-q^{40})-2q^4\psi(-q^{12})\psi(-q^{20})+2q^4\psi(q^2)\psi(q^{30}),\\
    &\hspace{6.5 cm} \text{[from (\ref{eq:11}) with $-q^4$]}\\
    2q^4\psi(-q^2)\psi(-q^{30})&=4q^4\psi(q^{12})\psi(q^{20})-2q^4\psi(q^2)\psi(q^{30}),\\
    4q^{16}\psi(q^8)\psi(q^{120})&=2q^4\psi(q^{12})\psi(q^{20})-2q^4\psi(-q^{12})\psi(-q^{20}).
\end{align*}
Putting these in (\ref{eq:15}) yields after rearrangement
\begin{align*}
    \varphi(q)\varphi(q^{15})-\varphi(q^6)\varphi(q^{10})+4q^4\psi(q^2)\psi(q^{30})-4q^4\psi(q^{12})\psi(q^{20})=2q\psi(q^3)\psi(q^5).
\end{align*}

Replacing this in (\ref{eq:9}), we have
\begin{align}
    \sum_{n=0}^\infty u_1(16n+6)q^n=16q\psi(q^3)\psi(q^3)\psi(q^5).\label{eq:14}
\end{align}
By (\ref{eq:4}) and (\ref{eq:14}), we get
\begin{align*}
    16\sum_{n=0}^\infty T(1,6,15;2n)q^{n+1}=\sum_{n=0}^\infty u_1(16n+6)q^n.
\end{align*}
Comparing the coefficients of $q^{n+1}$, we can conclude that
\begin{align*}
    16T(1,6,15;2n)=3N(1,6,15;16n+22)-N(1,6,15;4(16n+22)).
\end{align*}

\noindent\textbf{Part(ii)}: We will only prove Theorem \ref{thm4} (ii) for $(a,b,c)=(2,3,5)$. The other cases can be proved similarly using Lemma \ref{lemma2} and Lemma \ref{lemma1}.

By (\ref{eq:2}), we have
\begin{align*}
    \sum_{n=0}^\infty T(2,3,5;n)q^n=\psi(q^2)\psi(q^3)\psi(q^5).
\end{align*}
In view of (\ref{eq:8}), we get
\begin{align*}
    \sum_{n=0}^\infty T(2,3,5;2n+1)q^{2n+1}&=\psi(q^2)\mathcal{O}(\psi(q^3)\psi(q^5))\\
    &=\psi(q^2)(q^3\psi(q^2)\psi(q^{30})).
\end{align*}
Dividing by $q$ and substituting $q$ by $q^{1/2}$ yields after rearrangement
\begin{align}
    \sum_{n=0}^\infty T(2,3,5;2n+1)q^n=q\psi(q)\psi(q)\psi(q^{15}).\label{eq:22}
\end{align}
On the other hand, by (\ref{eq:1}), we have
\begin{align}
    \sum_{n=0}^\infty N(2,3,5;n)q^n=\varphi(q^2)\varphi(q^3)\varphi(q^5).\label{eq:23}
\end{align}
In view of (\ref{eq:5}), we get
\begin{align*}
    \sum_{n=0}^\infty N(2,3,5;n)q^n=(\varphi(q^8)+2q^2\psi(q^{16}))(\varphi(q^{12})+2q^3\psi(q^{24}))(\varphi(q^{20})+2q^5\psi(q^{40})).
\end{align*}
Extracting the terms involving $q^{4n}$ and substituting $q$ by $q^{1/4}$ yields after rearrangement
\begin{align}
    \sum_{n=0}^\infty N(2,3,5;4n)q^n=\varphi(q^2)\varphi(q^3)\varphi(q^5)+4q^2\varphi(q^2)\psi(q^6)\psi(q^{10}).\label{eq:24}
\end{align}
We define
\begin{align*}
    u_2(n):=3N(2,3,5;n)-N(2,3,5;4n).
\end{align*}
By (\ref{eq:23}) and (\ref{eq:24}), we have
\begin{align*}
    \sum_{n=0}^\infty u_2(n)q^n=2\varphi(q^2)\varphi(q^3)\varphi(q^5)-4q^2\varphi(q^2)\psi(q^6)\psi(q^{10}).
\end{align*}
In view of (\ref{eq:5}), we get
\begin{align*}
    \sum_{n=0}^\infty u_2(2n)q^{2n}&=2\varphi(q^2)\mathcal{E}(\varphi(q^3)\varphi(q^5))-4q^2\varphi(q^2)\psi(q^6)\psi(q^{10})\\
    &=2\varphi(q^2)(\varphi(q^{12})\varphi(q^{20})+4q^8\psi(q^{24})\psi(q^{40}))-4q^2\varphi(q^2)\psi(q^6)\psi(q^{10}).
\end{align*}
Substituting $q$ by $q^{1/2}$ yields after rearrangement
\begin{align*}
    \sum_{n=0}^\infty u_2(2n)q^n=2\varphi(q)\varphi(q^6)\varphi(q^{10})-4q\varphi(q)\psi(q^3)\psi(q^5)+8q^4\varphi(q)\psi(q^{12})\psi(q^{20}).
\end{align*}
In view of (\ref{eq:5}) and (\ref{eq:8}), we get
\begin{align*}
    \sum_{n=0}^\infty &u_2(4n+2)q^{2n+1}\\
    &=2\varphi(q^6)\varphi(q^{10})\mathcal{O}(\varphi(q))-4q\mathcal{E}(\varphi(q))\mathcal{E}(\psi(q^3)\psi(q^5))-4q\mathcal{O}(\varphi(q))\mathcal{O}(\psi(q^3)\psi(q^5))\\
    &\quad+8q^4\psi(q^{12})\psi(q^{20})\mathcal{O}(\varphi(q))\\
    &=2\varphi(q^6)\varphi(q^{10})(2q\psi(q^8))-4q(\varphi(q^4))(\psi(q^8)\varphi(q^{60})+q^{14}\varphi(q^4)\psi(q^{120}))\\
    &\quad-4q(2q\psi(q^8))(q^3\psi(q^2)\psi(q^{30}))+8q^4\psi(q^{12})\psi(q^{20})(2q\psi(q^8)).
\end{align*}
Dividing by $q$ and substituting $q$ by $q^{1/2}$ yields after rearrangement
\begin{align*}
    \sum_{n=0}^\infty u_2(4n+2)q^n&=-4\varphi(q^2)\psi(q^4)\varphi(q^{30})+4\varphi(q^3)\psi(q^4)\varphi(q^5)-8q^2\psi(q)\psi(q^4)\psi(q^{15})\\
    &\quad+16q^2\psi(q^4)\psi(q^6)\psi(q^{10})-4q^7\varphi(q^2)\varphi(q^2)\psi(q^{60}).
\end{align*}
In view of (\ref{eq:5}) and (\ref{eq:3}), we get
\begin{align*}
    \sum_{n=0}^\infty &u_2(8n+2)q^{2n}\\
    &=-4\varphi(q^2)\psi(q^4)\varphi(q^{30})+4\psi(q^4)\mathcal{E}(\varphi(q^3)\varphi(q^5))-8q^2\psi(q^4)\mathcal{E}(\psi(q)\psi(q^{15}))\\
    &\quad+16q^2\psi(q^4)\psi(q^6)\psi(q^{10})\\
    &=-4\varphi(q^2)\psi(q^4)\varphi(q^{30})+4\psi(q^4)(\varphi(q^{12})\varphi(q^{20})+4q^8\psi(q^{24})\psi(q^{40}))\\
    &\quad-8q^2\psi(q^4)(\psi(q^6)\psi(q^{10}))+16q^2\psi(q^4)\psi(q^6)\psi(q^{10}).
\end{align*}
Substituting $q$ by $q^{1/2}$ yields after rearrangement
\begin{align*}
    \sum_{n=0}^\infty u_2(8n+2)q^n&=-4\varphi(q)\psi(q^2)\varphi(q^{15})+4\psi(q^2)\varphi(q^6)\varphi(q^{10})+8q\psi(q^2)\psi(q^3)\psi(q^5)\\
    &\quad+16q^4\psi(q^2)\psi(q^{12})\psi(q^{20}).
\end{align*}
In view of (\ref{eq:5}) and (\ref{eq:8}), we get
\begin{align*}
    \sum_{n=0}^\infty &u_2(16n+2)q^{2n}\\
    &=-4\psi(q^2)\mathcal{E}(\varphi(q)\varphi(q^{15}))+4\psi(q^2)\varphi(q^6)\varphi(q^{10})+8q\psi(q^2)\mathcal{O}(\psi(q^3)\psi(q^5))\\
    &\quad+16q^4\psi(q^2)\psi(q^{12})\psi(q^{20})\\
    &=-4\psi(q^2)(\varphi(q^4)\varphi(q^{60})+4q^{16}\psi(q^8)\psi(q^{120}))+4\psi(q^2)\varphi(q^6)\varphi(q^{10})\\
    &\quad+8q\psi(q^2)(q^3\psi(q^2)\psi(q^{30}))+16q^4\psi(q^2)\psi(q^{12})\psi(q^{20}).
\end{align*}
Substituting $q$ by $q^{1/2}$ yields after rearrangement
\begin{align}
\begin{aligned}
    \sum_{n=0}^\infty &u_2(16n+2)q^n\\
    &=8q^2\psi(q)\psi(q)\psi(q^{15})\\
    &\quad+4\psi(q)\Big(-\varphi(q^2)\varphi(q^{30})+\varphi(q^3)\varphi(q^5)+4q^2\psi(q^6)\psi(q^{10})-4q^8\psi(q^4)\psi(q^{60})\Big).\label{eq:25}
\end{aligned}
\end{align}

In view of (\ref{eq:5}), we get
\begin{align*}
    &-\varphi(q^2)\varphi(q^{30})+\varphi(q^3)\varphi(q^5)+4q^2\psi(q^6)\psi(q^{10})-4q^8\psi(q^4)\psi(q^{60})\\
    &=-\varphi(q^2)\varphi(q^{30})\\
    &\quad+\Big(\varphi(q^{12})\varphi(q^{20})+2q^3\varphi(q^{20})\psi(q^{24})+2q^5\varphi(q^{12})\psi(q^{40})+4q^8\psi(q^{24})\psi(q^{40})\Big)\\
    &\quad+\Big(2q^2\psi(q^6)\psi(q^{10})+2q^2\psi(q^6)\psi(q^{10})\Big)-4q^8\psi(q^4)\psi(q^{60})\\
    &=2q^2\Big(\psi(q^6)\psi(q^{10})+q\varphi(q^{20})\psi(q^{24})+q^3\varphi(q^{12})\psi(q^{40})\Big)-\varphi(q^2)\varphi(q^{30})\\
    &\quad+\varphi(q^{12})\varphi(q^{20})+2q^2\psi(q^6)\psi(q^{10})-4q^8\psi(q^4)\psi(q^{60})+4q^8\psi(q^{24})\psi(q^{40})
\end{align*}
Now by (\ref{eq:3}), we have
\begin{align}
\begin{aligned}
    &-\varphi(q^2)\varphi(q^{30})+\varphi(q^3)\varphi(q^5)+4q^2\psi(q^6)\psi(q^{10})-4q^8\psi(q^4)\psi(q^{60})\\
    &=2q^2\psi(q)\psi(q^{15})-\varphi(q^2)\varphi(q^{30})+\varphi(q^{12})\varphi(q^{20})+2q^2\psi(q^6)\psi(q^{10})\\
    &\quad-4q^8\psi(q^4)\psi(q^{60})+4q^8\psi(q^{24})\psi(q^{40}).\label{eq:26}
\end{aligned}
\end{align}
Implementing (\ref{eq:11}), (\ref{eq:12}), (\ref{eq:13}), and (\ref{eq:10}) appropriately, we get
\begin{align*}
    \varphi(q^2)\varphi(q^{30})&=\varphi(-q^{12})\varphi(-q^{20})+2q^2\psi(q^6)\psi(q^{10}),\\
    \varphi(q^{12})\varphi(q^{20})&=\varphi(-q^8)\varphi(-q^{120})+2q^8\psi(q^4)\psi(q^{60})\\
    &=\varphi(-q^{12})\varphi(-q^{20})-2q^8\psi(-q^4)\psi(-q^{60})+2q^8\psi(q^4)\psi(q^{60}),\\
    &\hspace{6.5 cm}\text{[from (\ref{eq:12}) with $-q^4$]}\\
    4q^8\psi(q^{24})\psi(q^{40})&=2q^8\psi(q^4)\psi(q^{60})+2q^8\psi(-q^4)\psi(-q^{60}).
\end{align*}
Putting these in (\ref{eq:26}) yields after rearrangement
\begin{align*}
    -\varphi(q^2)\varphi(q^{30})+\varphi(q^3)\varphi(q^5)+4q^2\psi(q^6)\psi(q^{10})-4q^8\psi(q^4)\psi(q^{60})=2q^2\psi(q)\psi(q^{15}).
\end{align*}

Replacing this in (\ref{eq:25}), we have
\begin{align}
    \sum_{n=0}^\infty u_2(16n+2)q^n=16q^2\psi(q)\psi(q)\psi(q^{15}).\label{eq:27}
\end{align}
By (\ref{eq:22}) and (\ref{eq:27}), we get
\begin{align*}
    16\sum_{n=0}^\infty T(2,3,5;2n+1)q^{n+1}=\sum_{n=0}^\infty u_2(16n+2)q^n.
\end{align*}
Comparing the coefficients of $q^{n+1}$, we can conclude that
\begin{align*}
    16T(2,3,5;2n+1)=3N(2,3,5;16n+18)-N(2,3,5:4(16n+18)).
\end{align*}

\section{Acknowledgements}
The authors would like to thank IISER, TVM for providing the excellent working conditions.
The authors would like to thank T. Kathiravan for helpful comments and suggestions.

\end{document}